\documentclass[12pt,reqno]{article}
\usepackage[usenames]{color}
\usepackage{amssymb}
\usepackage{amsmath}
\usepackage{amsfonts}
\usepackage[T1]{fontenc}
\usepackage{enumerate}

\newtheorem{theorem}{Theorem}

\newtheorem{corollary}[theorem]{Corollary}

\allowdisplaybreaks
\begin{document}

\begin{center}
\vskip1cm

{\LARGE \textbf{Elementary function representations for the moments of the
Meyer-K\"{o}nig and Zeller operators}}

\vspace{2cm}

{\large Ulrich Abel}\\[3mm]
\textit{Technische Hochschule Mittelhessen}\\[0pt]
\textit{Fachbereich MND}\\[0pt]
\textit{Wilhelm-Leuschner-Stra\ss e 13, 61169 Friedberg, }\\[0pt]
\textit{Germany}\\[0pt]
\href{mailto:Ulrich.Abel@mnd.thm.de}{\texttt{Ulrich.Abel@mnd.thm.de}}
\end{center}

\vspace{2cm}

{\large \textbf{Abstract.}}

\bigskip

We calculate the moments of the Meyer-K\"{o}nig and Zeller operators in
terms of elementary functions and polylogarithms.

\bigskip



\vspace{2cm}

\section{Introduction}

The operators of Meyer-K\"{o}nig and Zeller \cite{Meyer-K-Zeller-1960} in
the slight modification of Cheney and Sharma \cite{Cheney-Sharma-1964},
defined by 
\begin{equation}
\left( M_{n}f\right) \left( x\right) =\left( 1-x\right) ^{n+1}\sum_{\nu
=0}^{\infty }{\binom{\nu +n}{\nu }}f\left( \frac{\nu }{\nu +n}\right) x^{\nu
}\text{ }\qquad \left( 0\leq x<1\right) ,  \label{Def-MKZ-univariate}
\end{equation}%
(also called Bernstein power series, or briefly MKZ\ operators) were the
object of several investigations in approximation theory.

Let $e_{r}\left( x\right) =x^{r}$ denote the monomials. The moments of the
MKZ operators are given by$\ \left( M_{n}e_{r}\right) \left( x\right) $.
They are starting point for several studies of their approximation
properties. In 1984, J. A. H. Alkemade \cite{Alkemade-JAT-1984} found a
representation of the second moment 
\begin{equation}
\left( M_{n}e_{2}\right) \left( x\right) =x^{2}+\frac{x\left( 1-x\right) ^{2}%
}{n+1}\text{ }_{2}^{\text{ }}F_{1}\left( 1,2;n+2;x\right)  \label{Alkemade}
\end{equation}%
in terms of hypergeometric functions. The way of deriving this expression is
based upon a differential equation which is satisfied by the right-hand side
of $\left( \ref{Def-MKZ-univariate}\right) $. Alkemade wrote:
\textquotedblleft For the MKZ operators, an explicit expression for $\left(
M_{n}e_{2}\right) \left( x\right) $ does not yet occur in the
literature.\textquotedblright\ \cite[Page 262]{Alkemade-JAT-1984}. Gavrea
and Ivan \cite[Theorem 2.2, Eq. (2.3)]{Gavrea-Ivan-MJOM-2018} observed that
the second moment is an elementary function and found an explicit
representation of $\left( M_{n}e_{2}\right) \left( x\right) $ in terms of
the logarithm and rational functions.

Entering $_{2}F_{1}\left( 1,2;n+2;x\right) $ into a computer algebra software 
(Wolfram Mathematica
Version 9), for small values
of the integer variable $n$ (up to $n=9$), one obtains explicit expressions
of the form 
\begin{equation*}
\sum_{k=1}^{n}\frac{a_{k}}{x^{k}}+\log \left( 1-x\right) \sum_{k=2}^{n+1}%
\frac{b_{k}}{x^{k}}.
\end{equation*}%
The command \texttt{Simplify} reveals the second sum to be a multiple of $%
\left( 1-x\right) ^{n+1}x^{-n-1}$. This is by no means surprising: 
\begin{equation*}
\text{ }_{2}F_{1}\left( 1,2;n+2;x\right) =\left( n+1\right)
!\sum_{k=0}^{\infty }\frac{\left( k+1\right) !}{\left( n+k+1\right) !}x^{k}.
\end{equation*}%
Now $\sum_{k=0}^{\infty }\frac{\left( k+1\right) !}{\left( n+k+1\right) !}%
x^{k}$ is a $n$-th antiderivative of $\sum_{k=0}^{\infty }x^{k+1}=\frac{1}{%
1-x}-1$. Therefore, it can be expected that $_{2}F_{1}\left(
1,2;n+2;x\right) $ has a representation in terms of logarithmic functions $%
\log \left( 1-x\right) $ and rational functions.

The purpose of this note is to provide a short and direct proof of this
fact. It follows by a more general representation for arbitrary MKZ moments
in terms of polylogarithmic functions. The proof uses only the definition $%
\left( \ref{Def-MKZ-univariate}\right) $ of the MKZ operators. We mention
that an integral representation for $\left( M_{n}e_{r}\right) \left(
x\right) $ can be found in \cite{Abel-MKZ moments-JAT-1995}. \ It is useful
when studying the asymptotic properties of the MKZ\ operators $\left(
M_{n}f\right) \left( x\right) $ for smooth functions $f$ as $n$ tends to
infinity (see \cite{Abel-JMAA-1997}).

\section{Results}

For $s=0,1,2,\ldots $, put $g_{s}\left( x\right) =\left( 1-x\right) ^{s}$.
Obviously, $M_{n}g_{0}=e_{0}$. For $s,n\in \mathbb{N}$, we obtain 
\begin{eqnarray*}
\left( M_{n}g_{s}\right) \left( x\right)  &=&\left( 1-x\right)
^{n+1}\sum_{\nu =0}^{\infty }{\binom{\nu +n}{\nu }}\left( \frac{n}{\nu +n}%
\right) ^{s}x^{\nu } \\
&=&\frac{n^{s}}{n!}\left( 1-x\right) ^{n+1}\sum_{\nu =0}^{\infty }\left( \nu
+n\right) \cdots \left( \nu +2\right) \left( \nu +1\right) \frac{x^{\nu }}{%
\left( \nu +n\right) ^{s}} \\
&=&\frac{n^{s}}{n!}\left( 1-x\right) ^{n+1}\left( \sum_{\nu =0}^{\infty }%
\frac{x^{\nu +n}}{\left( \nu +n\right) ^{s}}\right) ^{\left( n\right) } \\
&=&\frac{n^{s}}{n!}\left( 1-x\right) ^{n+1}\text{Li}_{s}^{\left( n\right)
}\left( x\right) ,
\end{eqnarray*}%
where 
\begin{equation*}
\text{Li}_{s}\left( x\right) =\sum_{\nu =1}^{\infty }\frac{x^{\nu }}{\nu ^{s}%
}\text{ \qquad }\left( \left\vert x\right\vert <1\right)
\end{equation*}%
is the polylogarithm. Using $e_{r}=\sum_{s=0}^{r}\left( -1\right) ^{s}\binom{%
r}{s}g_{s}$ we obtain the following result.

\begin{theorem}
For $r,n\in \mathbb{N}$ and $\left\vert x\right\vert <1$, the moments of the
MKZ operators have the representation 
\begin{equation*}
\left( M_{n}e_{r}\right) \left( x\right) =1+\frac{1}{n!}\left( 1-x\right)
^{n+1}\sum_{s=1}^{r}\left( -1\right) ^{s}\binom{r}{s}n^{s}\text{Li}%
_{s}^{\left( n\right) }\left( x\right) .
\end{equation*}
\end{theorem}

Applying the well-known (elementary) relation $\left( d/dx\right) $Li$%
_{s}\left( x\right) =x^{-1}$Li$_{s-1}\left( x\right) $, for $s\geq 2$, one
can reduce the terms Li$_{s}^{\left( n\right) }\left( x\right) $ to certain
derivatives of Li$_{1}\left( x\right) $. Note that Li$_{1}\left( x\right)
=-\log \left( 1-x\right) $ is an elementary function. As an immediate
consequence all moments of the MKZ operators can be represented in terms of
polylogarithms and rational functions. In particular, we have, for $n\in 
\mathbb{N}$, 
\begin{eqnarray*}
\text{Li}_{1}^{\left( n\right) }\left( x\right)  &=&\frac{\left( n-1\right) !%
}{\left( 1-x\right) ^{n}}, \\
\text{Li}_{2}^{\left( n\right) }\left( x\right)  &=&\left( x^{-1}\text{Li}%
_{1}\left( x\right) \right) ^{\left( n-1\right) }=\left( -1\right) ^{n}\frac{%
\left( n-1\right) !}{x^{n}}\log \left( 1-x\right)  \\
&&+\sum_{k=1}^{n-1}\left( -1\right) ^{n-1-k}\binom{n-1}{k}\frac{\left(
n-k-1\right) !}{x^{n-k}}\frac{\left( k-1\right) !}{\left( 1-x\right) ^{k}} \\
&=&\left( -1\right) ^{n}\frac{\left( n-1\right) !}{x^{n}}\log \left(
1-x\right) +\frac{\left( n-1\right) !}{x^{n}}\sum_{k=1}^{n-1}\left(
-1\right) ^{n-1-k}\frac{1}{k}\left( \frac{x}{1-x}\right) ^{k}.
\end{eqnarray*}%
In the special case $r=1$ we recover the well-known fact that MKZ operators
preserve linear functions:\ 
\begin{equation*}
\left( M_{n}e_{1}\right) \left( x\right) =1+\frac{1}{n!}\left( 1-x\right)
^{n+1}\text{Li}_{1}^{\left( n\right) }\left( x\right) =1-\left( 1-x\right)
=x.
\end{equation*}%
In the special case $r=2$ we obtain 
\begin{equation*}
\left( M_{n}e_{2}\right) \left( x\right) =1+\frac{1}{n!}\left( 1-x\right)
^{n+1}\left( -rn\text{Li}_{1}^{\left( n\right) }\left( x\right) +n^{2}\text{%
Li}_{2}^{\left( n\right) }\left( x\right) \right) 
\end{equation*}%
A short calculation leads to the following formula which is comparable to 
\cite[Theorem 2.2, Eq. (2.3)]{Gavrea-Ivan-MJOM-2018}.

\begin{corollary}
For $n\in \mathbb{N}$ and $\left\vert x\right\vert <1$, the second moment of
the MKZ operators has the representation 
\begin{eqnarray*}
\left( M_{n}e_{2}\right) \left( x\right) &=&\left( -1\right) ^{n}\frac{%
n\left( 1-x\right) ^{n+1}}{x^{n}}\log \left( 1-x\right) \\
&&+2x-1+n\left( 1-x\right) \sum_{k=1}^{n-1}\left( -1\right) ^{n-1-k}\frac{1}{%
k}\left( \frac{1-x}{x}\right) ^{n-k}.
\end{eqnarray*}
\end{corollary}

\strut

\thispagestyle{empty}


\begin{thebibliography}{9}
\bibitem{Abel-MKZ moments-JAT-1995} Ulrich Abel, \newblock The moments for
the Meyer-K\"{o}nig and Zeller operators, \newblock J. Approx. Theory 
\textbf{82} (1995), 352--361.

\bibitem{Abel-JMAA-1997} Ulrich Abel, The complete asymptotic expansion for
the Meyer-K\"{o}nig and Zeller operators, J. Math. Anal. Appl. \textbf{208}
(1997), 109--119.

\bibitem{Alkemade-JAT-1984} J. A. H.~Alkemade, \newblock The second moment
for the Meyer-K\"{o}nig and Zeller operators, \newblock J. Approx. Theory 
\textbf{40} (1984), 261--273.

\bibitem{Cheney-Sharma-1964} E. W. Cheney and A. Sharma, Bernstein power
series, Canad. J. Math. \textbf{16} (1964), 241--253.

\bibitem{Gavrea-Ivan-MJOM-2018} Ioan Gavrea and Mircea Ivan, \newblock An
elementary function representation of the second-order moment of the Meyer-K%
\"{o}nig and Zeller operators, \newblock Mediterr. J. Math. \textbf{15}
(2018), No. 1, Paper No. 20, 8 p.

\bibitem{Meyer-K-Zeller-1960} W. Meyer-K\"{o}nig and K. Zeller,
Bernsteinsche Potenzreihen, Studia Math. \textbf{19} (1960), 89--94.
\end{thebibliography}
\end{document}